\documentclass[12pt,leqno]{amsart}
\def\ds{\displaystyle}

\newtheorem{theorem}{Theorem}
\def\CC{{\mathbb C}}
\def\dist{\operatorname{dist}}
\def\supp{\operatorname{supp}}
\def\Ree{\operatorname{Re}}
\title[]
{Behavior of the Bergman kernel and metric near convex boundary points}
\author{Nikolai Nikolov and Peter Pflug }

\address
{Institute of Mathematics and Informatics\\
Bulgarian Academy of Sciences\\
1113 Sofia, Bulgaria}
\email{nik@math.bas.bg}
\address{Carl von Ossietzky Universit\"at Oldenburg\\
Fachbereich Mathematik\\ Postfach 2503\\
D-26111 Oldenburg,Germany}
\email{pflug@mathematik.uni-oldenburg.de}

\begin{document}


\begin{abstract}
The boundary behavior of the Bergman metric near a convex boundary
point $z_0$ of a pseudoconvex domain $D\subset\CC^n$ is studied;
it turns out that the Bergman metric at points $z\in D$ in
direction of a fixed vector $X_0\in\CC^n$ tends to infinite, when
$z$ is approaching $z_0$, if and only if the boundary of $D$ does
not contain any analytic disc through $z_0$ in direction of $X_0$.
\end{abstract}

\maketitle

For a domain $D\subset\CC^n$ we denote by $L_h^2(D)$ the Hilbert
space of all holomorphic functions $f$ that are square-integrable
and by $||f||_D$ the $L_2$-norm of $f$. Let $K_D(z)$ be the
restriction on the diagonal of the Bergman kernel function of $D$.
It is well-known (cf. \cite{JP}) that $$K_D(z)=\sup\{|f(z)|^2:f\in
L_h^2(D),\;||f||_D\le1\}.$$ If $K_D(z)>0$ for some point $z\in D$,
then the Bergman metric $B_D(z;X),\\ X\in\CC^n$, is well-defined
and can be given by the equality
$$B_D(z;X)=\frac{M_D(z;X)}{\sqrt{K_D(z)}},$$ where
$$M_D(z;X)=\sup\{|f'(z)X|:f\in L_h^2(D),\,||f||_D=1,\;f(z)=0\}.$$

We say that a boundary point $z_0$ of a domain $D\subset\CC^n$ is
{\it convex} if there is a neighborhood $U$ of this point such
that $D\cap U$ is convex.

In \cite{Her}, Herbort proved the following

\begin{theorem} Let $z_0$ be a convex boundary point of a bounded pseudoconvex
domain $D\subset\CC^n$ whose boundary contains no nontrivial germ
of an analytic curve near $z_0$. Then $$\lim_{z\to
z_0}B_D(z;X)=\infty$$ for any $X\in\CC^n\setminus\{0\}$.
\end{theorem}

Herbort's proof is mainly based on Ohsawa's
$\bar\partial$-technique. The main purpose of this note is to
generalize Theorem 1 using more elementary methods.

For a convex boundary point $z_0$ of a domain $D\subset\CC^n$ we
denote by $L(z_0)$ the set of all $X\in\CC^n$ for which there
exists a number $\varepsilon_X>0$ such that $z_0+\lambda
X\in\partial D$ for all complex numbers $\lambda$,
$|\lambda|\le\varepsilon_X$. Note that $L(z_0)$ is a complex
linear space.

Then our result is the following one.

\begin{theorem} Let $z_0$ be a convex boundary point of a bounded pseudoconvex
domain $D\subset\CC^n$ and let $X\in\CC^n$. Then

{\rm(a)} $\ds\liminf_{z\to z_0}K_D(z)\dist^2(z,\partial D)\in
(0,\infty]$;

{\rm(b)} $\ds\lim_{z\to z_0}B_D(z;X)=\infty$ if and only if
$X\not\in L(z_0)$. Moreover, this limit is locally uniform in
$X\not\in L(z_0)$;

{\rm(c)} If $L(z_0)=\{0\}$, then {\rm(a)} and {\rm(b)} are still
true without the assumption that $D$ is bounded.

\end{theorem}

{\bf Proof of Theorem 2.} To prove (a) and (b) we will use the
following localization theorem for the Bergman kernel and metric
\cite{DFH}.

\begin{theorem} Let $D\subset\CC^n$ be a bounded pseudoconvex domain and let
$V\subset\subset U$ be open neighborhoods of a point
$z_0\in\partial D$ Then there exists a constant $\tilde C\ge 1$
such that $$\tilde C^{-1}K_{D\cap U}(z)\le K_D(z)\le K_{D\cap
U}(z),$$ $$\tilde C^{-1}B_{D\cap U}(z;X)\le B_D(z;X)\le\tilde C
B_{D\cap U}(z;X)$$ for any $z\in D\cap V$ and any $X\in\CC^n$.
{\rm (}Here $K_{D\cap U}(z)$ and $B_{D\cap U}(z;\cdot)$ denote the Bergman
kernel and metric of the connected component of $D\cap U$ that
contains $z$.{\rm )}
\end{theorem}

So, we may assume that $D$ is convex.

To prove part (a) of Theorem 2, for any $z\in D$ we choose a point
$\tilde z\in\partial D$ such that $||z-\tilde z||=\dist(z,\partial
D)$. We denote by $l$ the complex line through $z$ and $\tilde z$.
By the Oshawa-Takegoshi extension theorem for $L^2$--holomorphic
functions \cite{OT}, it follows that there exists a constant
$C_1>0$ only depending on the diameter of $D$ (not on $l$) such
that $$K_D(z)\ge C_1 K_{D\cap l}(z).\leqno{(1)}$$ Since ${D\cap
l}$ is convex, it is contained in an open half-plane $\Pi$ of the
$l$-plane with $\tilde z\in\partial\Pi$. Then $$K_{D\cap l}(z)\ge
K_\Pi(z)=\frac{1}{4\pi \dist^2(z,\partial\Pi)}.\leqno(2)$$ Now,
part (a) of Theorem 2 follows from the inequalities (1), (2) and
the fact that  $\dist(z,\partial\Pi)\le||z-\tilde
z||=\dist(z,\partial D)$.

To prove part (b) of Theorem 2, we denote by $N(z_0)$ the complex
affine space through $z_0$ that is orthogonal to $L(z_0)$. Set
$E(z_0)=D\cap N(z_0)$. Note that $E(z_0)$ is a nonempty convex
set. So, part (b) of Theorem 2 will be a consequence of the
following

\begin{theorem} Let $z_0$ be a boundary point of a bounded convex domain
$D\subset\CC^n$. Then:

{\rm(i)} $\ds\lim_{z\to z_0}B_D(z;X)=\infty$ locally uniformly in
$X\not\in L(z_0)$;

{\rm(ii)} for any compact set $K\subset\subset E(z_0)$ there
exists a constant $C>0$ such that $$B_D(z;X)\le C||X||,\quad z\in
K^0,\;X\in L(z_0),$$ where $K^0:=\{z_0+t z:\ z\in K,\ 0<t\le 1\}$
is the cone generated by $K$.
\end{theorem}

{\bf Proof of Theorem 4.} To prove (i) we will use the well-known
fact that the Carath\'eodory metric $C_D(z;X)$ of $D$ does not
exceed $B_D(z;X)$. On the other hand, we have the following simple
geometric inequality \cite{BP}: $$C_D(z;X)\ge\frac{1}{2d(z;X)},$$
where $d(z;X)$ denotes the distance from $z$ to the boundary of
$D$ in the $X$-direction, i.e. $d(z;X):=\sup\{r:z+\lambda X\in
D,\;\lambda\in\CC,\;|\lambda|<r\}$. So, if we assume that (i) does
not hold, then we may find a number $a>0$ and sequences
$D\supset(z_j)_j,\; z_j\to z_0,\ \CC^n\supset(X_j)_j,\;X_j\to
X\not\in L(z_0)$ such that $\ds B_D(z_j;X_j)\le\frac{1}{2a}$.
Hence $d(z_j;X_j)\ge a$ which implies that for $|\lambda|\le a$
the points $z_0+\lambda X$ belong to $\bar D$ and, in view of
convexity, they belong to $\partial D$. This means that $X\in
L(z_0)$ -- contradiction.

To prove part (ii) of Theorem 4, we may assume that $z_0=0$ and
$L:=L(0)=\{z\in\CC^n:z_1=\ldots=z_k=0\}$ for some $k<n$. Then
$N:=N(0)=\{z\in\CC^n:z_{k+1}=\ldots=z_n=0\}$. From now on we will
write any point $z\in\CC^n$ in the form $z=(z',z'')$,
$z'\in\CC^k$, $z''\in\CC^{n-k}$. Note that $L\in\partial D$ near
$0$, i.e. there exists an $c>0$ such that
$$\{0'\}\times\Delta''_c\in\partial D,\leqno{(3)}$$ where
$\Delta''_c\subset\CC^{n-k}$ is the polydisc with center at the
origin and radius $c$. Since $K\subset\subset E:=E(0)$ and since
$E$ is convex, there exists an $\alpha>1$ such that
$K\subset\subset E_\alpha$, where $E_\alpha:=\{z:\alpha z\in E\}$.
Note that $K^0\subset E_\alpha$. Using (3), the following equality
$$(z',z'')=\frac{1}{\alpha}(\alpha z',0'')+(1-\frac{1}{\alpha})
(0',(1-\frac{1}{\alpha})^{-1}z''),$$ and the convexity of $D$, it
follows that $$F_\alpha\times\Delta''_\varepsilon\subset
D,\leqno{(4)}$$ where $\ds\varepsilon:=c(1-\frac{1}{\alpha})$ and
where $F_\alpha$ is the projection of $E_\alpha$ in $\CC^k$ (we
can identify $E_\alpha$ with $F_\alpha$). For
$\delta:=c(\alpha-1)$ we get in the same way that $$\tilde
D:=D\cap(\CC^k\times\Delta''_\delta)\subset
F_{\frac{1}{\alpha}}\times\Delta''_\delta.\leqno{(5)}$$

Now, let $(z,X)\in K^0\times L$. Note that $z=(z',0'')$ and
$X=(0',X'')$. Then, using (4) and the product properties of the
Bergman kernel and metric, we have $$M_D(z;X)\le
M_{F_\alpha\times\Delta''_\varepsilon}(z;X)=\leqno{(6)}$$
$$=M_{\Delta''_\varepsilon}(0'';X'')\sqrt{K_{F_\alpha}(z')}\le
C_1||X||\sqrt{K_{F_\alpha}(z')}$$ for some constant $C_1>0$. On
the other side, since
$K^0\subset\subset\CC^k\times\Delta''_\delta$, in virtue of
Theorem 3 there exists a constant $\tilde C\ge 1$ such that
$$K_D(z)\ge\tilde C^{-1}K_{\tilde D}(z).$$ Moreover, in view of
(5), we have $$K_{\tilde D}(z)\ge
K_{F_{\frac{1}{\alpha}}}(z')K_{\Delta''_\delta}(0'')$$ and hence
$$K_D(z)\ge (C_2)^2 K_{F_{\frac{1}{\alpha}}}(z')\leqno{(7)}$$ for
some constant $C_2>0$. Now, by (6) and (7), it follows that
$$B_D(z;X)=\frac{M_D(z;X)}{\sqrt{K_D(z)}}\le\frac{C_1}{C_2}||X||
\sqrt{\frac{K_{F_\alpha}(z')}{K_{F_{\frac{1}{\alpha}}}(z')}}.\leqno{(8)}$$
Note that $z'\to\alpha^{-2}z'$ is a biholomorphic mapping from
$F_{\frac{1}{\alpha}}$ onto $F_\alpha$ and, therefore,
$$K_{F_{\frac{1}{\alpha}}}(z')=\alpha^{-4k}K_{F_\alpha}(\alpha^{-2}z').\leqno{(9)}$$
In view of (8) and (9), in order to finish (ii) we have to find a
constant $C_3>0$ such that $$K_{F_\alpha}(z')\le C_3
K_{F_\alpha}(\alpha^{-2}z')\leqno{(10)}$$ for any $z'\in H^0$ with
$H^0:=\{t z':\ z'\in H,\ 0<t\le 1\}$, where $H$ is the projection
of $K$ into $\CC^k$ (we can identify $K$ with $H$).

To do this, note first that $\gamma:=\dist(H,\partial F_\alpha)>0$
since $K\subset\subset E_\alpha$. Fix $\tau\in(0,1]$ and $z'\in
H^0$, and denote by $T_{\tau,z'}$ the translation that maps the
origin in the point $\tau z'$. It is easy to check that
$$T_{\tau,z'}(\bar F_\alpha\cap B_\gamma)\subset
F_\alpha,\leqno{(11)}$$ where $B_\gamma$ is the ball in $\CC^k$
with center at the origin and radius $\gamma$. To prove (10), we
will consider the following two cases:

Case I. $z'\in H^0\setminus B_{\frac{\gamma}{2}}\subset\subset
F_\alpha$: Then
$$K_{F_\alpha}(z')\le\frac{m_1}{m_2}K_{F_\alpha}(\alpha^{-2}z'),\leqno{(12)}$$
where $\ds m_1:=\sup_{H^0\setminus
B_\frac{\gamma}{2}}K_{F_\alpha}$ and $m_2:=\inf K_{F_\alpha}$.

Case II. $z'\in H^0\cap B_\frac{\gamma}{2}$: By Theorem 3 there
exists a constant $\tilde C_3\ge 1$ such that $\tilde C_3
K_{F_\alpha}\ge K_{F_\alpha\cap B_\gamma}$ on $\ds F_\alpha\cap
B_{\frac{\gamma}{2}}$. In particular, $$\tilde
C_3K_{F_\alpha}(\alpha^{-2}z')\ge K_{F_\alpha\cap
B_\gamma}(\alpha^{-2}z').\leqno{(13)}$$ On the other side, by (11)
with data $T:=T_{1-\alpha^{-2},z'}$ it follows that
$$K_{F_\alpha\cap B_\gamma}(\alpha^{-2}z')=K_{T(F_\alpha\cap
B_\gamma)}(z')\ge K_{F_\alpha}(z').\leqno(14)$$ Now, (12), (13),
and (14) imply that (10) holds for $\ds
C_3:=\max\{\frac{m_1}{m_2},\tilde C_3\}$. This completes the
proofs of Theorem 4 and part (b) of Theorem 2.

\qed

{\bf Remark.} The approximation (5) of the domain
$D\cap(\CC^k\times\Delta''_\delta)$ by the domain $\ds
E_{\frac{1}{\alpha}}\times\Delta''_\delta$ can be replaced by
using the Oshawa-Takegoshi theorem \cite{OT} with the data $D$ and
$N$.

Finally, part (c) of Theorem 2 will be a consequence of the
following two theorems.

\begin{theorem} \cite{Nik} Let $D\subset\CC^n$. be a pseudoconvex domain and let
$U$ be an open neighborhood of a local (holomorphic) peak point
$z_0\in\partial D$. Then $$\lim_{z\to z_0}\frac{K_D(z)}{K_{D\cap
U}(z)}=1$$ and $$\lim_{z\to z_0}\frac{B_D(z;X)}{B_{D\cap
U}(z;X)}=1$$ locally uniformly $X\in\CC^n\setminus\{0\}$.
\end{theorem}

\begin{theorem} Let $z_0$ be a boundary point of a bounded convex domain
$D\subset\CC^n$. Then the following conditions are equivalent:

1. $z_0$ is a (holomorphic) peak point;

2. $z_0$ is the unique analytic curve in $\bar D$ containing
$z_0$;

3. $L(z_0)=\{0\}$.
\end{theorem}

Note that the only nontrivial implication is
$(3)\Longrightarrow(1)$ it is contained in \cite{Sib}. Now, part
(c) of Theorem 2 is a consequence of this implication, Theorem 5,
and part (b) of Theorem 2.

{\bf Proof of Theorem 6.} The implication $(2)\Longrightarrow(3)$ is trivial.

The implication $(1)\Longrightarrow(2)$ easily follows by the
maximum principle and the fact that there are a neighborhood $U$
of $z_0$ and a vector $X\in\CC^n$ such that $(\bar D\cap
U)+(0,1]X\subset D$ (cf. (11)).

Denote by $A^0(D)$ the algebra of holomorphic functions on $D$
which are continuous on $\bar D$. Now, following \cite{Sib} we
shall prove the implication $(3)\Longrightarrow(1)$; namely, (3)
implies that  $z_0$ is a peak point with respect to $A^0(D)$. This
is equivalent to the fact (cf. \cite{Gam}) that the point mass at
$z_0$ is the unique element of the set $A(z_0)$ of all
representing measures for $z_0$ with respect to $A^0(D)$, i.e.
$\supp\ \mu=\{z_0\}$ for any $\mu\in A(z_0)$.

Let $\mu\in A(z_0)$. Since $D$ is convex, we may assume that
$z_0=0$ and $D\subset\{z\in\CC^n:\Ree(z_1)<0\}$. Note that if $a$
is a positive number such that $\ds a\inf_{z\in D}\Ree(z_1)>-1$ (D
is bounded), then the function $f_1(z)=\exp(z_1+az_1^2)$ belongs
to $A^0(D)$ and $|f_1(z)|<1$ for $z\in\bar D\setminus\{z_1=0\}$.
This easily implies (cf. \cite{Gam}) that $\supp\ \mu\subset
D_1:=\partial D\cap\{z_1=0\}$. Since $L(0)=0$, the origin is a
boundary point of the compact convex set $D_1$. As above, we may
assume that $D_1\subset\{z\in\CC^n:\Ree(z_2)\le 0\}$ ($z_2$ is
independent of $z_1$) and then construct a function $f_2\in
A^0(D)$ such that $|f_2(z)|<1$ for $z\in D_1\setminus\{z_2=0\}$.
This implies that $\supp\ \mu\subset D_1\cap\{z_2=0\}$. Repeating
this argument we conclude that $\supp\ \mu=\{0\}$, which completes
the proofs of Theorems 6 and 2.

\qed

{\bf Acknowledgments.} This paper was written during the stay of
the first author at the University of Oldenburg supported by a
grant from the DAAD. He likes to thank both institutions, the DAAD
and the Mathematical Department of the University of Oldenburg.

\end{document}